\documentclass[psamsfonts]{amsart}

\usepackage{amsmath,amssymb}

\newtheorem{prop}[subsection]{Proposition}
\newtheorem{cor}[subsection]{Corollary}

\newtheorem{rmk}[subsection]{Remark}
      {\begin{rmk}}
      {\end{rmk}}

\newcommand{\cal}{\mathcal}
\newcommand{\be}{\begin}
\newcommand{\en}{\end}

\newcommand{\down}{\textsf{$\downarrow^B$}}
\newcommand{\mathbfr}{\mathrm}
\newcommand{\Pow}{\textsf{Pow}}

\newcommand{\CZF}{{\mathbfr{CZF}}}

\newcommand{\wREA}{{\mathbfr{wREA}}}

\newcommand{\comment}[1]{}

\makeatletter
\def\blfootnote{\xdef\@thefnmark{}\@footnotetext}
\makeatother

\begin{document}

\title{On Tarski's fixed point theorem}
\author{Giovanni Curi}
\address{}
\curraddr{}
\email{giovanni.curi@email.it}
\thanks{}


\date{}

\dedicatory{To Orsola}


\begin{abstract}
A concept of \emph{abstract inductive definition} on a complete lattice is formulated and studied.
As an application, a constructive  version of Tarski's fixed point theorem is obtained.
\end{abstract}

\maketitle


\section*{Introduction}
The fixed point theorem referred to in this paper is the one asserting that every monotone mapping on a complete lattice $L$ has a least fixed point. The  proof, due to A. Tarski, of this result, is
a simple and  most significant example
of a proof that can be carried out on the base of intuitionistic logic (e.g. in the intuitionistic set theory IZF, or in topos logic), and that yet
is widely regarded as essentially non-constructive.

The reason for this fact is that Tarski's construction of the fixed point is highly impredicative: if $f:L\to L$ is a monotone map, its least fixed point is given by $\bigwedge P$, with $P\equiv \{x\in L \mid f(x)\leq x\}$. Impredicativity here is found in the fact that the fixed point, call it $p$, appears in its own construction ($p$ belongs to  $P$), and, indirectly, in the fact that the complete lattice $L$ (and, as a consequence, the  collection $P$ over which the infimum is taken) is assumed to form a set,  an assumption that seems only reasonable in an intuitionistic setting in the presence of strong impredicative principles (cf. Section \ref{uniformobjects} below).

In concrete applications (e.g. in computer science and numerical analysis) the monotone operator $f$ is often also continuous, in particular it preserves suprema of non-empty chains; in this situation, the least fixed point can be constructed taking the supremum of the ascending chain $\bot, f(\bot), f(f(\bot)),...,$ given by the set of finite iterations of $f$ on the least element $\bot$.
This procedure has also been extended to the general (non-continuous) case, but using transfinite iterations of $f$ along the classical ordinals, and strongly non-constructive arguments.

Several alternative, more constructive, proofs of Tarski's theorem have  been proposed,  probably the two most satisfactory so far being those obtained in intuitionistic contexts in  \cite{JM,Taylor}. These proofs, however, still presuppose  the existence of powerobjects (use the powerset axiom), and/or make use of fully impredicative comprehension principles.

In this paper I present a constructive  predicative version of Tarski's fixed point theorem. Working in the context of Myhill-Aczel's constructive set theories\footnote{These theories are constructive  subtheories of classical ZF set theory.}, I formulate a concept of \emph{abstract inductive definition} on a complete lattice, and use it to obtain a generalization  to complete lattices of Aczel's  theory of inductive definitions on a set.
Every abstract inductive definition gives rise to an inductively defined  `generalized element' of the lattice. Using this fact, I  derive first a proof of Tarski's  theorem  in a basic system for constructive set theory extended by the impredicative full Separation scheme. This proof may  be regarded as an improvement on the mentioned results in \cite{JM,Taylor}. Then, under assumptions on the lattice and the monotone map that are always satisfied in fully impredicative (classical or intuitionistic) systems, a constructive  predicative proof of Tarski's  theorem is  obtained.

Before discussing (respectively in sections \ref{aid} and \ref{Tarski}) abstract inductive definitions and their application to the fixed point theorem, I formulate in Section \ref{uniformobjects} a notion of \emph{uniform class}, by analogy with the concept of uniform object studied in the context of the effective topos \cite{Hy}. This notion  is used proof-theoretically to show that in the systems for constructive set theory we work with, certain standard partially ordered structures, as complete  lattices or directed-complete partial orders, must be defined as having a proper class of elements.

\section{Constructive set theory}\label{CZF}
\setcounter{subparagraph}{0}
I shall be working in the setting of Myhill-Aczel's constructive set theories. The basic system in this context  is the choice-free  \emph{Constructive Zermelo-Fraenkel Set Theory} (CZF), due to Aczel.
This system is often extended by principles, such as the Regular Extension Axiom, that ensure that certain inductively defined classes are sets, and with choice principles, as Dependent Choice (DC) or the stronger Presentation Axiom (PA).
Note that CZF, extended or not by these principles, is a subtheory of classical set theory ZFC. In  contrast to ZF and its intuitionistic version IZF, CZF does not have the impredicative unrestricted Separation scheme and the Powerset axiom.
Here I  provide the basic information that make the paper self-contained; the reader may consult \cite{AR} for a thorough introduction to the subject.

The language of CZF is the same as that of Zermelo-Fraenkel Set Theory, ZF, with $\in$ as the only non-logical symbol. Beside the rules and axioms of a standard calculus for intuitionistic predicate logic with equality, CZF has the following axioms and axiom schemes:

\begin{enumerate}
  \item Extensionality: $\forall a\forall b (\forall y (y \in a \leftrightarrow y \in b) \rightarrow a = b)$.
  \item Pair: $\forall a\forall b\exists x\forall y (y \in x\leftrightarrow y = a \vee y = b).$
  \item Union: $\forall a\exists  x\forall y (y \in x \leftrightarrow  (\exists z \in a)(y \in z))$.
  \item Restricted Separation scheme: \smallskip

\qquad \qquad $\forall a\exists  x\forall y (y \in x \leftrightarrow  y\in a\wedge \phi(y)),$
 \smallskip

 \noindent for $\phi$ a restricted formula. A formula $\phi$ is \emph{restricted} if the quantifiers that occur in it are of the form $\forall x\in b$,  $\exists x\in c$.

  \item Subset Collection scheme:
  \begin{eqnarray*}&&\forall a\forall b\exists c\forall u ((\forall x \in a)(\exists y \in b)\phi(x, y, u)\; \rightarrow\\
  &&  (\exists d \in c)((\forall x \in a)(\exists y \in d)\phi (x, y, u)
\wedge (\forall y \in d)(\exists x \in a)\phi (x, y, u))).\end{eqnarray*}

  \item Strong Collection scheme:
  \begin{eqnarray*}&&\forall a ((\forall x \in a)\exists y \phi(x, y)\;\rightarrow\\
      && \exists b ((\forall x \in  a)(\exists y \in b)\phi (x, y) \wedge (\forall y\in b)(\exists x \in a)\phi(x, y))).
   \end{eqnarray*}

  \item Infinity: $\exists a (\exists x\in a  \wedge  (\forall x\in a)(\exists y\in a)  x\in y)$.

  \item Set Induction scheme: $\forall a ((\forall x \in a)\phi(x)\rightarrow \phi(a))\rightarrow \forall a\phi(a)$.
\end{enumerate}

\noindent
We shall denote by $\CZF^-$ the system obtained from $\CZF$ by leaving out the Subset Collection scheme.
Subset Collection is perhaps the most unusual of the CZF axioms and schemes; for  this paper it suffices to note that using it one proves that the class $b^a$ of functions  from a set $a$ to a set $b$ is a set, i.e., the  Exponentiation Axiom.

Friedman's Intuitionistic Zermelo-Fraenkel set theory based on collection,  IZF, has the same theorems as $\CZF$ extended by the unrestricted Separation Scheme and the Powerset Axiom. Moreover, the theory obtained from $\CZF$, or from IZF, by adding the Law of Excluded Middle has the same theorems as ZF.

As in classical set theory, one makes use in this context of class notation and terminology \cite{AR}. For example, given any set or class $X,$ one has the class $\Pow(X)=\{x\mid x\subseteq X\}$ of subsets of $X$. A major role in constructive set theory is played by inductive definitions.  Inductively defined classes and sets  are natural objects in a constructive context  \cite{A86}; further, they partially compensate the lack of the impredicative powerset construction.

An \emph{inductive definition}  is any class  $\Phi$ of pairs.
A class $A$ is \emph{$\Phi-$closed} if:
\smallskip

\centerline{$(a,X)\in \Phi$, and  $X\subseteq A$ implies $a\in A$.}
\smallskip

\noindent
The following  theorem is called \emph{the class inductive definition theorem} \cite{AR}.

\be{thm}[$\CZF^-$]\label{cidt}
Given any class  $\Phi$ of ordered pairs, there exists a least $\Phi-$closed class $I(\Phi)$, the class inductively defined by $\Phi$.
\en{thm}

\noindent An \emph{inductive definition on a set} is an inductive definition $\Phi$ that is a subclass of the cartesian product  $S\times \Pow(S)$, for $S$ a set.

Even when $\Phi$ is a set, $I(\Phi)$ need not be a set in CZF. For this reason, CZF is often extended with the  \emph{Regular Extension Axiom}, REA.
\vspace{0.3mm}

\qquad REA: every set is a subset of a regular set.
\vspace{0.3mm}

\noindent A set $c$ is \emph{regular} if it is transitive, inhabited, and for any $u\in c$ and any set $R\subseteq u \times c$,  if $(\forall x \in u)(\exists y) \langle x, y\rangle \in R$, then there is a set $v\in  c$ such that
\begin{eqnarray}\label{rea}(\forall x \in u)(\exists y \in v)((x, y) \in R) &\wedge& (\forall y \in v)(\exists x\in u)((x, y)  \in R).\end{eqnarray}
$c$ is said to be {\em weakly regular} if in the above definition of regularity the second conjunct in (\ref{rea}) is omitted.
The weak regular extension axiom, wREA, is the statement that every set is the subset of a weakly regular set. In $\CZF + \wREA$, the following  theorem can be proved.

\be{thm}[$\mathbfr{CZF + wREA}$]\label{sidt}
If  $\Phi$ is a set, then $I(\Phi)$ is a set.
\en{thm}

\noindent
More generally, Theorem \ref{sidt} holds for inductive definitions that are \emph{bounded}  \cite{AR}.

CZF extended by (w)REA  has been given a fully constructive justification \cite{A86},
and has shown to provide a reasonably adequate system
for the development of intuitionistic predicative mathematics (e.g. \cite{A,AR,CuExistenceSC}). Note that CZF +(w)REA is a subsystem of classical set theory ZFC \cite{RathjenLubarsky2003}.

\section{The uniform objects in constructive set theory}
\label{uniformobjects}
The notion of uniform object was introduced in connection with the effective topos  (cf. \cite{Hy,Ro}).
In the context of  set theory, we shall say that a class $A$ is \emph{uniform} if, for every set $a$, and every formula $\phi(x,y)$,\smallskip

\qquad \qquad $(\forall x\in A)(\exists y\in a)\phi(x,y)\to(\exists y\in a)(\forall x\in A)\phi(x,y).$
\smallskip

\noindent
Note that if $A$ is a uniform class, and $f:A\to a$, with $a$ any set, is a function, then $f$ must be constant.

Every singleton set (terminal object) is uniform. If a set is not a singleton, it  is not uniform, as follows logically by the definition of uniformity: define a set $c$ \emph{trivial} if it is a singleton, i.e. if $(\exists x\in c)(\forall y\in c)(x=y)$; c is \emph{non-trivial} if it is not trivial.

\be{prop}\label{uniformset}
Every uniform set is trivial.
\en{prop}

\be{proof} Assume $A$ is a uniform class, and that there is  a set  $a$ with $a=A$. Then, by
uniformity of $A$, choosing
$\phi(x,y)\equiv x=y$  one gets $(\exists x\in a)(\forall y\in a)(x=y)$.
\en{proof}

\noindent Note that the empty set is non-trivial, and therefore  non-uniform.

Proposition \ref{uniformset} holds (in particular) for CZF and every extension of it.
In fact, CZF,  as well as any extension $T$ of CZF consistent with the law of excluded middle LEM, cannot prove that a non-trivial class $A$ is uniform. Indeed, if $T$ proves $A$ uniform, $T$ + LEM proves $A$ trivial: since the empty set is not uniform, there must be $a\in A$. Then, if $b\in A$ is such that $b\neq a$, a non-constant function $f$ from $A$ to $\{0,1\}$ can be defined ($f$ sends $a$ to $0$,  and every element of $A$ different from $a$ to $1$), against the uniformity of $A$. Thus, $b=a$ for all $b\in A$, so that $A$ is trivial.

However, in extensions of CZF that are \emph{not} consistent with LEM, one does find non-trivial uniform classes (even in these extensions, by Proposition \ref{uniformset}, one will not find non-trivial uniform \emph{sets}).
Recall that CZF is consistent with the conjunction of the following two non-classical principles.
\smallskip

\noindent \emph{Uniformity principle}, \textbf{UP}: for every formula $\phi$,

\begin{itemize}
  \item[] $(\forall x)(\exists y\in \mathbb{N})\phi(x,y)\to(\exists y\in \mathbb{N})(\forall x)\phi(x,y).$
\end{itemize}

\noindent  \emph{Every set is subcountable}, \textbf{SC}:

\begin{itemize}
  \item []  $(\forall x)(\exists U\in \Pow(\mathbb{N}))(\exists f)f:U\twoheadrightarrow x,$
\end{itemize}

\noindent where $f:U\twoheadrightarrow x$ indicates that $f$ is surjective from $U$ to $x$.

The uniformity principle is a generalization of the principle that any function from the class of all sets to $\mathbb{N}$ must be constant. The informal justification for this  principle is that the class of all sets forms such an `indiscernible'  totality that the only total functions from it to a set as $\mathbb{N}$ can be the constant ones (cf. \cite[pg. 234]{TvD} for more on this point).
SC asserts  that every set can be enumerated by a subsets of the natural numbers; the idea behind it is the intuitive one that sets have to be no larger than $\mathbb{N}$.\footnote{The reader may feel uncomfortable with these two principles. We hasten to say that they will only be exploited in this section to obtain underivability results.}

 Using these two principles one finds a first interesting uniform class.

\begin{prop}[CZF+UP+SC]\label{Vuniform}
The universal class $\mathrm{V}=\{x\mid x=x\}$ is uniform.
\end{prop}

\be{proof}  Let $a$ be a set and $\phi$ be a formula.  Assume that for all $x\in \mathrm{V}$ there exists $y\in a$ such that $\phi(x,y)$. By SC, there are $U\in  \Pow(\mathbb{N})$, and  $f:U\twoheadrightarrow a$. Thus,  for all $x\in \mathrm{V}$ there exists $n\in \mathbb{N}$ such that $[n\in U\ \&\ \phi(x,f(n))]$. One can then apply UP, and easily conclude.
\en{proof}

\noindent Uniformity of the universal class can also be regarded as a principle in itself,
the \emph{generalized uniformity principle} (GUP). In \cite{CuriPA, CuExistenceSC} we exploited this general form of the uniformity principle  for proving the independence  from constructive set theory of various (classically and) intuitionistically valid results. Models of CZF (and of several of its extensions, including the system CZF + Sep + REA + PA), that also validate the conjunction of UP and SC, have been described by various authors \cite{BergMoerdijk,Lubarsky2006,Rathjen94,Rathjen06,Streicher}.

In the following,  CZF$^\sharp$ denotes any possible extension of CZF that is simultaneously  consistent with SC and UP, and so with GUP.
\smallskip

\noindent As for the uniform objects in the effective topos (cf. \cite{Ro}),  we have:

\be{lem}\label{surjection}
For classes $X,Y, f$, if $X$ is uniform, and $f:X\twoheadrightarrow Y$ is an onto mapping, then $Y$ is uniform (and therefore a proper class if non-trivial).
\en{lem}

\be{proof} Easy calculations.
\en{proof}

\noindent Using this lemma, we may find new interesting uniform classes.

\be{cor}\label{PowUnif}
CZF + UP + SC proves that for every set $X$, $\Pow(X)$ is uniform. The same system plus Sep (full Separation) proves that $\Pow(X)$ is uniform for every class $X$.
\en{cor}

\be{proof} Assume $X$ is a set. Then the map $f:\mathrm{V}\to \Pow(X)$, defined by $f(y)=y\cap X$, is onto. With Sep, the same holds even if $X$ is a class, as by Sep one gets that $y\cap X$ is a set and is therefore in $\Pow(X)$. By \ref{surjection} and \ref{Vuniform} one concludes.
\en{proof}

One might be tempted to think that uniformity of an object is connected with its size, and that in particular, if an injection $X\hookrightarrow Y$ of a uniform class $X$ into $Y$ exists, then $Y$ is uniform. This is not so, as one may realize
considering the class $\Pow(\{0\})\cup \{*\}$, with $*$ an element not belonging to $\Pow(\{0\})$: clearly a non-constant function can  be defined on this class.

\smallskip

A partially ordered class, or \emph{poclass}, $(X,\leq)$ is a class $X$ together with a class-relation $\leq$ that is reflexive, transitive, and antisymmetric. $X$ is \emph{flat} if $(\forall x,y\in X)(x\leq y\to x=y)$.

A partially ordered class $(X,\leq)$ is a (large) \emph{$\bigvee$-semilattice} if every subset has a supremum; $(X,\leq)$ is a (large) \emph{$\bigwedge$-semilattice} if it has infima of arbitrary subsets
(note: a large $\bigvee$-semilattice need not be a $\bigwedge$-semilattice, nor conversely). A poclass  is \emph{directed complete}, or a (large) \emph{dcpo}, if it has joins of \emph{directed} subsets (a subset $U$ of a poclass $X$ is directed if it is inhabited, and whenever $x,y\in U$ there is $z\in U$ with $x,y\leq z$). A poclass  is \emph{conditionally  complete}, or a (large) \emph{bcpo} if every \emph{inhabited} and \emph{bounded} subset has a join. Finally, a partially ordered class is \emph{chain-complete}, or a (large) \emph{ccpo},  if every chain (i.e., totally ordered subset) has a join.

Our main application of the concept of uniform class is that, whenever non-flat, these structures cannot be carried by sets in CZF$^\sharp$. Note that if a poclass with a top or bottom element is flat, then it is trivial. So for  $\bigvee$-semilattices, $\bigwedge$-semilattices and chain-complete poclasses, to be non-flat it suffices to be non-trivial, as these structures always possess a bottom, or top, element.

Part of the following result is proved in \cite{CuExistenceSC}; the proof to be presented gives in particular a more uniform explanation of the results obtained there.

\be{thm}\label{largestructures}
The following partially ordered structures are carried by proper classes in CZF + SC + UP, and thus cannot be proved to be carried by sets in CZF, and its extensions CZF$^\sharp$:
\begin{enumerate}
  \item  Non-trivial $\bigvee$-semilattices and $\bigwedge$-semilattices.
  \item \label{dcpos} Non-flat dcpo's; in particular, non-trivial dcpo's with a top or a bottom element.

  \item \label{bcpos} Non-flat bcpo's; in particular, non-trivial bcpo's with a top or a bottom element.
  \item \label{ccpos} Non-trivial chain-complete partially ordered classes.
\end{enumerate}

\en{thm}

\be{proof}
1. The result for $\bigvee$-semilattices is proved  in \cite{CuExistenceSC}. Note also that a $\bigvee$-semilattice is  a chain-complete partially ordered class, and, when non-trivial, a non-flat dcpo and bcpo. However,
as we shall be concerned with $\bigvee$-semilattices in the rest of this paper, we give a direct proof of this fact: let $L$ be a  $\bigvee$-semilattice. Assume $L$ is a set. Under this assumption, we may apply the first part of Corollary \ref{PowUnif} to get that $\Pow(L)$ is uniform. Then, since $\bigvee: \Pow(L)\to L$ is onto, $L$ is uniform too, by Lemma \ref{surjection}. Thus, by Proposition \ref{uniformset}, $L$ must be trivial. The proof for $\bigwedge$-semilattices is similar.

2.
Let $D$ be a dcpo, and let $b\leq a$, for $a,b\in D$. Assume $D$ is a set. Then, by Restricted Separation, the class \smallskip

\qquad \qquad $U_{a,b}(y)\equiv \{x\in D\mid x=a\ \&\ \emptyset \in y\}\cup\{b\}$
\smallskip

\noindent  is a set, and is directed, for every $y\in \mathrm{V}$. One may therefore consider the class
\smallskip

\qquad \qquad $K_{a,b}=\{U_{a,b}(y)\mid y\in \mathrm{V}\}.$
\smallskip

\noindent   $K_{a,b}$ is a uniform class (the map sending $y\in \mathrm{V}$ to $U_{a,b}(y)$ is onto). Then, again by the assumption that $D$ is a set, $\bigvee:K_{a,b}\to D$ has to be constant. Therefore
$\bigvee U_{a,b}(\emptyset)=b=a=\bigvee U_{a,b}(\{\emptyset\})$,
and this for every $b\leq a$, so that $D$ is flat.

3.  The claim is proved as for the previous case, since $U_{a,b}(y)$ is also bounded.

4. Let $C$ be non-trivial and chain-complete. $C$ has bottom $\bot=\bigvee \emptyset$. Assume $C$ is a set.  Considering, for $z\in C$, the uniform class of chains \smallskip

\qquad \qquad $K_z=\{W_z(y)\mid y\in \mathrm{V}\},\ \textrm{ with}\ W_z(y)\equiv \{x\in C\mid x=z \ \&\ \emptyset \in y\},$
\smallskip

\noindent   one may conclude reasoning similarly to the previous cases. Alternatively, one could observe that $U_{a,b}(y)$ is a chain for every $a,b$, and that a non-trivial ccpo is non-flat.
\en{proof}

The reader may be puzzled by this result: isn't for instance the  Boolean algebra $\{0,1\}$ in particular a complete lattice, and carried by a set? While $\{0,1\}$ is indeed a Boolean algebra (and is of course carried by a set), it is not intuitionistically complete, not even on the base of topos logic, or IZF (completeness of $\{0,1\}$ implies the weak law of excluded middle, e.g. \cite{FS79}).
More generally, one cannot expect finite lattices to be complete in an  intuitionistic setting.

Note that any set with the flat order $=$ is a dcpo and a bcpo (so that, requiring in items \ref{dcpos}, \ref{bcpos}, the weaker condition that the dcpo's, bcpo's are non-empty (or even inhabited) and non-trivial is not enough), and that a set with the same order and a bottom (or a top) element added is not a dcpo constructively.

In topos theory, and sometimes in constructive mathematics, one  considers the \emph{MacNeille reals} $R_m$ (also called \emph{extended reals} by Troelstra), see e.g. \cite{Elephant2}.
$R_m$ is a conditionally complete class, and,  in a topos, or in IZF, is a set.

\be{cor} The system CZF,  as well as  every extension  CZF$^\sharp$, cannot prove that
$R_m$ is  a set. The same holds for every conditionally complete extension of $\mathbb{Q}$.
\en{cor}

\noindent CZF is instead strong enough to prove that the Dedekind reals (that are not conditionally complete, intuitionistically)  form a set \cite{AR}.

\be{rmk}\rm\label{Jsuniform}
In CZF+UP+SC+Sep, one may improve on Theorem \ref{largestructures} (1). One can show indeed that any $\bigvee$-semilattice (or $\bigwedge$-semilattice) $L$ is actually a uniform class: by the second part of Corollary \ref{PowUnif}, in this context $\Pow(L)$ is uniform even with $L$ a class; then, as $\bigvee: \Pow(L)\to L$ is onto, $L$ is uniform by Lemma \ref{surjection}
(this is more generally the case for the carrier of every structure $(A, f,...)$ with at least one mapping $f:\Pow(A)\to A$ that is onto).
The same holds without Sep for so-called set-generated $\bigvee$-semilattices (see the next section): if $B$ is a generating set for $L$, the map $\bigvee: \Pow(B)\to L$ is again onto. By contrast, dcpo's, bcpo's, as well as their class of directed/bounded subsets, may not be uniform, even when they are non-flat:  the poclass $(\Pow(\{0\})\cup \{*\}, \leq)$, where $\leq$ is defined extending the inclusion relation on $\Pow(\{0\})$ by letting $*\leq *$, is a non-flat dcpo and bcpo.
\en{rmk}

By Theorem \ref{largestructures},   partially ordered structures of several types, that  are in particular structures of the types contemplated by Theorem \ref{largestructures}, fail to be carried by sets in CZF,  and in every extension  CZF$^\sharp$: these include various types of domains considered in the denotational semantics of programming languages, frames (locales), preframes,  continuous lattices (appropriately redefined), etc.

An appropriate way to deal with these structures in a constructive predicative system  is to regard them as \emph{generated} by sets; we shall in particular work with set-generated $\bigvee$-semilattices \cite{AR}, defined in the next section. In the presence of the powerset construction, set-generated $\bigvee$-semilattices are the same as standard $\bigvee$-semilattices, and similarly for the other structures.

\section{Abstract inductive definitions}
\label{aid}

In this and the next section abstract inductive definitions are introduced and used to extend the theory of inductive definitions on a set   (cf. \cite{A82,A86,AR}, or Section \ref{CZF}) to  $\bigvee$-semilattices. In this section we shall be working in CZF$^-$.

By Theorem \ref{largestructures}, in constructive set theory it is of no use to consider $\bigvee$-semilattices carried by sets.
The standard counterpart of the classical notion of $\bigvee$-semilattice is in this context given by the concept of \emph{set-generated $\bigvee$-semilattice} \cite{AR}.
A (large) $\bigvee$-semilattice $L$ is said to be set-generated if it has a generating set $B$, i.e. a subset $B$ of $L$ such that, for all  $x\in L$,
\begin{itemize}

\item [$i.$] $\down x\equiv \{b\in B\ \mid \ b\leq x\}$ is a set,

\item [$ii.$] $x=\bigvee \down x$.

\end{itemize}

\noindent The powerclass $\Pow(X)$ of a set $X$, ordered by inclusion, is a prototypical example of a set-generated $\bigvee$-semilattice, with generating set $B=\{\{x\}:x\in X\}$. Other examples are: the class of fixed points of a closure operator on $\Pow(X)$, for $X$ a set;
the Dedekind-MacNeille completion of a partially ordered set; the lattice of ideals on a distributive lattice carried by a set; the lattices of open subsets of any topological space with an explicitly given set-indexed base, and  the complete Boolean algebra  of $\neg\neg -$ stable elements of any such lattice \cite{CuriPA}. The notion of set-generated $\bigvee$-semilattice makes in fact possible to deal with a wide class of complete lattices in predicative systems, cf. \cite{AR,A} for more details. Note  that a set-generated $\bigvee$-semilattice is also a complete lattice.

\smallskip

Let $L$ be a set-generated $\bigvee$-semilattice $L$ with generating set $B$. An \emph{abstract inductive definition on $L$} (in the following often just an inductive definition) is any class of ordered pairs \smallskip

\qquad \qquad $\Phi\subseteq B\times L.$
\smallskip

\noindent To define what it means for a subclass of the generating set $B$ to be $\Phi-$closed we need the following notion.
A subclass $Y\subseteq B$ will be called \emph{\texttt{$c_L$}-closed} if, for every subset $U$ of $Y$, the set $\down \bigvee U$ is contained in $Y$;
i.e., $Y$ is \texttt{$c_L$}-closed if
\smallskip

\qquad \qquad $\bigcup_{U\in \Pow(Y)}\down \bigvee U = Y.$
\smallskip

\noindent
A $c_L$-closed class can be thought of as denoting a generalized element of $L$.\footnote{When $L$ is a complete Heyting algebra, such generalized elements also arise as generalized truth values in (predicatively defined) cHa-models for constructive set theory, cf. \cite{G1}.} If $Y$ is a set, $Y$ is \texttt{$c_L$}-closed iff  $Y=\down\bigvee Y$.

\noindent A  class $Y\subseteq B$ will be said \emph{$\Phi-$closed} if it is \texttt{$c_L$}-closed and if,  whenever $(b,a)\in \Phi$, \smallskip

\qquad \qquad $\down a\subseteq Y\ \Longrightarrow\ b\in Y.$\smallskip

\noindent
We shall denote by ${\mathcal I}(\phi)$  the least $\Phi-$closed class, if it exists.

\smallskip

Given an abstract inductive definition $\Phi$ on $L$, and an element $a$ in $L$,  the class
\smallskip

\qquad \qquad $\{b\in B \ \mid \ (\exists a')\ (b,a')\in \Phi \ \& \ a'\leq a\}$\smallskip

\noindent may not be a set in general. If for, every $a\in L$, this class is a set we say that $\Phi$ is \emph{local}. A  local abstract inductive definition $\Phi$  determines a mapping $\Gamma_\Phi: L\to L$, given by, for $a\in L$,
\smallskip

\qquad \qquad $\Gamma_\Phi (a)\equiv \bigvee\{b\in B \ \mid \ (\exists a')\ (b,a')\in \Phi \ \& \ a'\leq a\}.$\smallskip

\noindent If $a_1\leq a_2$, then $\Gamma_\Phi (a_1)\leq \Gamma_\Phi (a_2)$, i.e. $\Gamma_\Phi$ is \emph{monotone}.

Any monotone operator on $L$ can in fact be obtained in this way from a local abstract inductive definition.

\be{prop}\label{gamma=gammaPhi}
Let $\Gamma:L\to L$ be a monotone operator on $L$. Then, there is a local abstract inductive definition $\Phi_\Gamma$ such that, for every $a\in L$, $\Gamma(a)=\Gamma_{\Phi_\Gamma} (a)$.
\en{prop}

\be{proof} Define $\Phi_\Gamma\subseteq B\times L$ by \smallskip

\qquad \qquad $(b,a)\in \Phi_\Gamma\iff b\leq \Gamma(a).$\smallskip

\noindent $\Phi_\Gamma$ is local as, for  $a\in L$, $\{b\in B \ \mid \ (\exists a')\ (b,a')\in \Phi_\Gamma \ \& \ a'\leq a\}=\{b\in B \ \mid \ (\exists a')\  b\leq \Gamma(a')\ \& \ a'\leq a\}.$ By monotonicity of $\Gamma$, this class is the same as $\{b\in B \ \mid  b\leq \Gamma(a)\}$, that is a set by the assumption that $L$ is set-generated. The join of this set therefore exists, and again as $B$ is a set of generators, is equal to $\Gamma(a)$.
\en{proof}

We presented the above simple proof in detail with the purpose of emphasizing the role of the assumption that $L$ is set-generated. The following proposition explains the relationship between  $\Phi-$closed \emph{sets} of a local inductive definition,  and (pre-) fixed points of the associated monotone operator.

\be{prop}\label{idprefix}
Given a local inductive definition $\Phi$ on a $\bigvee$-semilattice $L$ with generating set $B$,
a one-to-one correspondence exists between the  $\Phi-$closed subclasses of $B$ that are sets, and the elements $a$ of $L$ such that $\Gamma_\Phi(a)\leq  a$. Moreover, whenever the class ${\cal I}(\Phi)$ exists and is a set, $\Gamma_\Phi$ has a least fixed point.
\en{prop}

\be{proof} Assume $Y\subseteq B$ is $\Phi-$closed and that it is a set. Then, $\bigvee Y$ exists in $L$ and we have:  \smallskip

\qquad \qquad  $\Gamma_\Phi(\bigvee Y)= \bigvee\{b\in B \ \mid \ (\exists a)\ (b,a)\in \Phi \ \& \ a\leq \bigvee Y\}.$\smallskip

\noindent To conclude that  $\Gamma_\Phi(\bigvee Y)\leq \bigvee Y$, let $b\in B$ be such that there is $a\in L$ with $(b,a)\in \Phi$ and $a\leq \bigvee Y$. To show that $b\leq \bigvee Y$ it suffices to prove that $\down a\subseteq Y$, since $Y$ is $\Phi-$closed. But this follows, by $a\leq \bigvee Y$, from the assumption that $Y$ is a set and that it is \texttt{$c_L$}-closed  (take $U=Y$ in the definition of \texttt{$c_L$}-closed).

Conversely, to $a\in L$ such that $\Gamma_\Phi(a)\leq  a$, we associate the \texttt{$c_L$}-closed class $\down a$. As $L$ is set-generated, $\down a$ is a set; using the assumption that $\Gamma_\Phi(a)\leq  a$, one immediately sees that $\down a$ is also $\Phi$-closed.

Finally, one has $Y=\down \bigvee Y$, as $Y$ is \texttt{$c_L$}-closed, and $a=\bigvee \down a$, since $L$ is set-generated.

Now assume ${\cal I}(\Phi)$ is a set. As ${\cal I}(\Phi)$ is $\Phi-$closed, by what has just been shown, $\Gamma_\Phi(\bigvee {\cal I}(\Phi))\leq  \bigvee {\cal I}(\Phi)$. To prove the converse, note that by monotonicity of $\Gamma_\Phi$,
$\Gamma_\Phi(\Gamma_\Phi(\bigvee {\cal I}(\Phi)))\leq  \Gamma_\Phi(\bigvee {\cal I}(\Phi))$. Then, $\down \Gamma_\Phi(\bigvee {\cal I}(\Phi))$ is $\Phi-$closed, again by the proved correspondence. Then, as ${\cal I}(\Phi)$ is the least $\Phi-$closed class, ${\cal I}(\Phi)\subseteq \down \Gamma_\Phi(\bigvee {\cal I}(\Phi))$, so that $\bigvee {\cal I}(\Phi)\leq  \Gamma_\Phi(\bigvee {\cal I}(\Phi))$. Thus $\bigvee {\cal I}(\Phi)$ is a fixed point for $\Gamma_\Phi$. If $a\in L$ is another fixed point, then in particular $\Gamma_\Phi(a)\leq  a$. Therefore $\down a$ is $\Phi-$closed, and  ${\cal I}(\Phi)\subseteq \down a$, which gives $\bigvee {\cal I}(\Phi)\leq  a$.
\en{proof}

When the class ${\cal I}(\Phi)$ exists and is a set, we shall refer to its join $a=\bigvee {\cal I}(\Phi)\in L$  as to the element inductively defined by $\Phi$.

\section{Tarski's fixed point theorem}
\label{Tarski}

Recall by the preliminaries that, by results in \cite{A82, AR},  for any (standard) inductive definition $\Phi\subseteq X\times \Pow(X)$, for $X$ any set, the least $\Phi-$closed class $I(\Phi)$ exists in the system  CZF$^-$. In this section we prove that more generally the least (\texttt{$c_L$}-closed and) $\Phi-$closed class ${\cal I}(\Phi)$ exists for every abstract inductive definition $\Phi$ on a set-generated $\bigvee$-semilattice $L$.

As a corollary of this result we shall have that CZF$^-$ extended by the impredicative unrestricted Separation scheme proves Tarski's fixed point theorem.
This improves on previous intuitionistic proofs of the theorem, as those in \cite{JM} or \cite{Taylor}.

Albeit `more constructive' than those proofs, this result can hardly be considered satisfactory.
For the case of monotone operators on $\bigvee$-semilattices of the form $\Pow(X)$ for $X$ a set, it directly follows by the results in \cite{A82, A86, AR} that a restricted version of Tarski's theorem obtains:
if a monotone operator $\Gamma: \Pow(X)\to \Pow(X)$ may be obtained as $\Gamma_\Phi$ by a bounded inductive definition $\Phi$ (in particular by an inductive definition $\Phi$ that is a set), then the system CZF + wREA proves that a least fixed point exists, as it proves that $I(\Phi)$ is a set in this case.

Here I show that, again in CZF + wREA, a bounded abstract inductive definition  on a set-generated $\bigvee$-semilattice $L$ gives rise to a least $\Phi-$closed class ${\cal I}(\Phi)$ that is a set,  whenever $L$ satisfies the extra standard condition of being set-presented. Thus, if a monotone operator $\Gamma:L\to L$ may be obtained as $\Gamma_\Phi$ for $\Phi$ a bounded abstract inductive definition, then $\Gamma$ has a fixed point.
In a sense there is no restriction here: in (I)ZF say, every abstract inductive definition is a set (and is therefore bounded),  every complete lattice is set-presented, and every monotone operator is obtained by a bounded inductive definition.
\smallskip

We now prove that for every abstract inductive definition $\Phi$ on $L$ the least (\texttt{$c_L$}-closed and) $\Phi-$closed class ${\cal I}(\Phi)$ exists. The classical proof would construct ${\cal I}(\Phi)$ iterating a monotone operator associated with $\Phi$ along the class of ordinals. For the case of monotone operators on point-classes, Aczel \cite{A82} showed that one can replace the class of ordinals with the class of all sets (applications of transfinite induction are then replaced by applications of Set Induction). Here we follow the same approach.

The proof that the least  $\Phi-$closed class ${\cal I}(\Phi)$ exists is a generalization of the proof in \cite{A82, AR} that the least $\Phi-$closed class exists for every standard inductive definition $\Phi$, considered  for the case $\Phi\subseteq X\times \Pow(X)$, for $X$ a set. Due to the `pointless' nature of $\bigvee$-semilattices some  difficulties however arise; the generalization in particular will require additional applications of the Strong Collection scheme. To make  visible what is involved in the generalization, we shall follow the terminology and notation of \cite{A82, AR} as far as possible; we shall however use suggestively Greek letters for sets playing the role of ordinals.
\smallskip

Let $B$ be a generating set for $L$ and $\Phi$ an inductive definition on $L$. Given a class of ordered pairs $Y\subseteq \mathrm{V}\times B$, let, for $\alpha$ any set in the universal class $\mathrm{V}$:\smallskip

\begin{itemize}
\item[] $Y^\alpha\equiv\{y\in B\mid (\alpha, y)\in Y\}$;\smallskip

\item[] $Y^{\in \alpha}\equiv\{y\in B\mid \exists \beta \in \alpha\ (\beta, y)\in Y\}\equiv \bigcup_{\beta\in \alpha} Y^\beta$;\smallskip

\item[] $Y^\infty\equiv \{y\in B\mid \exists \alpha \in \mathrm{V}\ (\alpha, y)\in Y\}\equiv \bigcup_{\alpha\in \mathrm{V}} Y^\alpha$.
\end{itemize}\smallskip

We shall need to consider an extension of the operator $\Gamma_\Phi$ to the generalized elements of $L$. To this aim, we define an operator on subclasses of $B$ which has the properties of a closure operator. For $Y$ a subclass of $B$, let
\smallskip

\qquad  $\texttt{$c_L$} Y \equiv \{b\in B\mid  (\exists U\in \Pow(B)) U\subseteq Y\ \& \ b\in \down \bigvee U\}= \bigcup_{U\in \Pow(Y)} \down \bigvee U$\smallskip

\noindent (cf. the extension $J$ of a $j-$operator on a set-generated frame in \cite{G1}). In the case of standard inductive definitions this operator has no role, as it reduces to identity.

The following proposition gathers together properties of $c_L$ that we shall need in the following. Its proof is a simple exercise  (due to the fact it is about classes, rather than just sets, in some cases Strong Collection is needed, cf. also \cite{G1}).

\be{prop}\label{c-operator} Let $X,Y, X_i$, for $i$ in a set $I$, be subclasses of $B$. Then
\begin{enumerate}
  \item $Y$ is \texttt{$c_L$}-closed iff $Y = \texttt{$c_L$} Y$;
  \item if $Y$ is a set, $\texttt{$c_L$} Y =  \down \bigvee Y;$

  \item $X\subseteq Y$ implies $\texttt{$c_L$} X\subseteq \texttt{$c_L$}Y$;
  \item $Y\subseteq \texttt{$c_L$} Y$;

  \item $\texttt{$c_L$}\texttt{$c_L$}Y = \texttt{$c_L$} Y$, so that \texttt{$c_L$}Y is \texttt{$c_L$}-closed;
  \item $X\subseteq \texttt{$c_L$}Y$ implies $\texttt{$c_L$} X\subseteq \texttt{$c_L$}Y$;
  \item $\texttt{$c_L$} \bigcup_{i\in I} X_i = \texttt{$c_L$} \bigcup_{i\in I} \texttt{$c_L$} X_i$.

\end{enumerate}
\en{prop}

\noindent
For $Y\subseteq B$, define \smallskip

\qquad $\bar \Gamma_\Phi (Y) \equiv \texttt{$c_L$} \{b\in B \ \mid \ (\exists a)\ (b,a)\in \Phi \ \& \ \down a\subseteq Y\}.$\smallskip

\noindent $\bar \Gamma_\Phi$ is clearly a monotone operator on classes.
\smallskip

\smallskip

Let $\Phi$ be an abstract inductive definition on a $\bigvee$-semilattice $L$ set-generated by a set $B$.

\be{lem}\label{classJ}
A class $J\subseteq \mathrm{V}\times B$ of ordered pairs exists such that, for every set $\alpha\in \mathrm{V}$,

\qquad \qquad \qquad \qquad $\texttt{$c_L$} J^\alpha = \bar \Gamma_\Phi (\texttt{$c_L$} J^{\in \alpha}).$
\en{lem}

\be{proof}
A set $G\subseteq \mathrm{V}\times B$ will be called \emph{good} (or an approximation to $J$) if, for any set $\alpha$,
\smallskip

\qquad \qquad $G^\alpha\subseteq \bar \Gamma_\Phi (\texttt{$c_L$} G^{\in \alpha}).$\smallskip

\noindent Define then \smallskip

\qquad \qquad $J=\bigcup \{G \mid G \ \textrm{good }\}.$
\smallskip

\noindent To prove $\texttt{$c_L$} J^\alpha \subseteq \bar \Gamma_\Phi (\texttt{$c_L$} J^{\in \alpha})$, it is enough to show that
$J^\alpha \subseteq \bar \Gamma_\Phi (\texttt{$c_L$} J^{\in \alpha})$, since $\bar \Gamma_\Phi(Y)$ is \texttt{$c_L$}-closed for every $Y$. So let $b\in J^\alpha$. Then a good set $G$ exists such that $b\in
\bar \Gamma_\Phi (\texttt{$c_L$} G^{\in \alpha}).$ Since $G^{\in \alpha}\subseteq J^{\in \alpha}$, also $\texttt{$c_L$} G^{\in \alpha}\subseteq \texttt{$c_L$}J^{\in \alpha}$, whence, by monotonicity of $\bar \Gamma_\Phi$, $b\in \bar \Gamma_\Phi (\texttt{$c_L$} J^{\in \alpha})$.

\smallskip

\noindent To prove the converse, let $y\in \bar \Gamma_\Phi (\texttt{$c_L$} J^{\in \alpha}).$ There is then a set $U\in \Pow(B)$ with \smallskip

\qquad \qquad $U\subseteq \{b\in B \ \mid \ (\exists a)\ (b,a)\in \Phi \ \& \ \down a\subseteq \texttt{$c_L$} J^{\in \alpha}\},\ \textrm{and} \ y\leq \bigvee U.$\smallskip

\noindent We shall prove that $U\subseteq J^\alpha$, so that $y\in  \texttt{$c_L$} J^\alpha$.
For  $b\in U$ there is $a\in L$ such that $(b,a)\in \Phi$ and  $ \down a \subseteq \texttt{$c_L$} J^{\in \alpha}$. The last inclusion can be rewritten as
\smallskip

\qquad  $(\forall y'\in \down a)\ (\exists W\in \Pow(J^{\in \alpha}))\ y'\leq \bigvee W$.
\smallskip

\noindent By Strong Collection, a set $K$ of subsets of $J^{\in \alpha}$ then exists such that
\smallskip

\qquad  $(\forall y'\in \down a)\ (\exists W\in K)\ y'\leq \bigvee W$.
\smallskip

\noindent Thus,
 $(\forall y'\in \down a)\ y'\leq \bigvee \bigcup K$, and $\bigcup K\subseteq J^{\in \alpha}$.
\smallskip

\noindent The latter expression can be reformulated as

\smallskip

\qquad \qquad $(\forall d\in \bigcup K)(\exists x\in \alpha)\ d\in J^x$, i.e.,

\smallskip

   \qquad \qquad $(\forall d\in \bigcup K)(\exists G)\ G\ \textrm{good}\ \& \ d\in  G^{\in \alpha}$.
\smallskip

\noindent By Strong Collection again, we then get a set $Z$ of good sets such that
\smallskip

\qquad \qquad $(\forall d\in \bigcup K)(\exists G\in Z)\ d\in  G^{\in \alpha}$.

\smallskip

\noindent Thus, $\bigcup K\subseteq  (\bigcup Z)^{\in \alpha}$.
Let then
\smallskip

\qquad \qquad \qquad $G=\{(\alpha,b)\}\cup \bigcup Z.$
\smallskip

\noindent As $\bigcup Z$ is a union of good sets, it is a good set too. Moreover, $\down a\subseteq \texttt{$c_L$} G^{\in \alpha}$, indeed, as seen, $(\forall y'\in \down a)\ y'\leq \bigvee \bigcup K$, and $\bigcup K\subseteq (\bigcup Z)^{\in \alpha}\subseteq G^{\in \alpha}$. Therefore, $b\in \bar \Gamma_\Phi (\texttt{$c_L$} G^{\in \alpha})$. Thus, $G$ is a good set.

Now, since $\{(\alpha, b)\}\in G$, we  have $b\in J^\alpha$. As this is true for every $b\in U$, and  $y\leq \bigvee U$, we get $y\in  \texttt{$c_L$} J^\alpha$, as wished.

\en{proof}

\be{thm}\label{leastidclass}
Let $\Phi$ be an abstract inductive definition on a $\bigvee$-semilattice $L$ set-generated by a set $B$.
Then, the smallest $\Phi-$closed class ${\cal I}(\Phi)$ exists.

\en{thm}

\be{proof} Our goal is to prove that the class \smallskip

\qquad \qquad $\texttt{$c_L$} J^\infty= \texttt{$c_L$}\bigcup_{\alpha\in \mathrm{V}}J^\alpha$
\smallskip

\noindent
is the least $\Phi-$closed class. Let $(b,a)\in \Phi$, and $\down a\subseteq \texttt{$c_L$} J^\infty$.  Then, given $y\in \down a$ there is a set $U$ such that
\smallskip

\qquad \qquad
$U\in \Pow(J^\infty)\ \& \ y\leq \bigvee U$.
\smallskip

\noindent
So,

\qquad \qquad $(\forall z \in U)(\exists \beta) z \in J^\beta$.
\smallskip

\noindent
By Collection, a set $\alpha$ exists such that
\smallskip

\qquad \qquad $(\forall z \in U)(\exists \beta\in \alpha) z \in J^\beta$.
\smallskip

\noindent
Thus, $U\subseteq J^{\in \alpha}$, so that $y\in \texttt{$c_L$} J^{\in \alpha}$.
So we have shown
\smallskip

\qquad \qquad \qquad \qquad
$(\forall y\in \down a)(\exists \alpha) y\in \texttt{$c_L$} J^{\in \alpha}.$
\smallskip

\noindent
Applying again the Collection scheme, we then get a set $K$ such that
\smallskip

\qquad \qquad \qquad \qquad $(\forall y\in \down a)(\exists \alpha\in K) y\in \texttt{$c_L$} J^{\in \alpha}.$
\smallskip

\noindent Since $\alpha \in K$ implies $\alpha\subseteq \bigcup K$, using Proposition \ref{c-operator} we conclude $\down a\subseteq \texttt{$c_L$} J^{\in \bigcup K}$. By definition of $ \bar \Gamma_\Phi$, then $b\in  \bar \Gamma_\Phi(\texttt{$c_L$} J^{\in \bigcup K})$. By Lemma \ref{classJ}, $\bar \Gamma_\Phi(\texttt{$c_L$} J^{\in \bigcup K}) = \texttt{$c_L$} J^{\bigcup K}$, so that, as $J^{\bigcup K}\subseteq J^\infty$, we conclude that $\texttt{$c_L$} J^\infty$ is $\Phi-$closed.

It remains to show that $\texttt{$c_L$} J^\infty$ is the least $\Phi-$closed class. Let $I$ be a $\Phi-$closed class, $I\subseteq B$. It is enough to prove that $J^\alpha\subseteq I$ for every set $\alpha$, since then, as $I$ is assumed  \texttt{$c_L$}-closed, by Proposition \ref{c-operator} one has $\texttt{$c_L$} \bigcup_{\alpha \in \mathrm{V}} J^\alpha = \texttt{$c_L$} J^\infty \subseteq I$.  We prove by Set Induction that, for every set $\alpha$, $J^\alpha\subseteq I$: let  $\alpha$ be a set. The inductive hypothesis gives $(\forall \beta\in\alpha) J^\beta\subseteq I$, i.e $J^{\in \alpha}\subseteq I$. By Proposition \ref{c-operator}, then $\texttt{$c_L$} J^{\in \alpha}\subseteq \texttt{$c_L$} I = I$. By monotonicity of  $\bar \Gamma_\Phi$, $\bar \Gamma_\Phi(\texttt{$c_L$} J^{\in \alpha})\subseteq \bar \Gamma_\Phi (I)\subseteq I$, whence, by Lemma \ref{classJ}, $J^\alpha\subseteq \texttt{$c_L$}J^\alpha\subseteq I$.
\en{proof}

\be{rmk}\rm\label{abstractstandard}
As observed by an anonymous referee, we could have deduced the above result by associating with the abstract inductive definition $\Phi$ the standard inductive definition $\Phi^*=\{(b,U)\in B\times \Pow(B)\mid b\leq \bigvee U\}\cup \{(b,\down a) \mid (b,a)\in \Phi\}$, and applying \cite[Theorem 5.1]{AR} (i.e., Theorem \ref{cidt} of Section \ref{CZF}). We remark however that the present  proof, beside being  direct and self-contained,  provides a generalization of the \emph{proof} of \cite[Theorem 5.1]{AR}.
In particular, it shows how to define constructively the transfinite iterations of a monotone operator on a general set-generated $\bigvee$-semilattice.
\en{rmk}

\noindent It may be worth noting that, so far, we made no use of Exponentiation (let alone of the Subset Collection scheme), so that the above results in fact hold in CZF$^-$.

For the following corollary we assume the impredicative unrestricted Separation scheme.

\begin{cor}[Tarski's fixed point theorem in CZF$^-$+ Sep]\label{TarskiSeparation}
The system CZF$^-$ augmented with the  Separation Scheme proves that every monotone operator $\Gamma$ on a set-generated  $\bigvee$-semilattice $L$ has a least fixed point.
\end{cor}

\be{proof} Let $\Phi_\Gamma$ be the inductive definition on $L$ associated with $\Gamma$ (see Section \ref{aid}). By the previous theorem, the least  $\Phi_\Gamma-$closed class ${\cal I}(\Phi_\Gamma)\subseteq B$ exists. Since $B$ is a set, by Separation  ${\cal I}(\Phi_\Gamma)$ is a set too. Then, by Propositions \ref{gamma=gammaPhi} and \ref{idprefix}, $\bigvee {\cal I}(\Phi_\Gamma)$ is the least fixed point of $\Gamma$.
\en{proof}

\noindent An elegant intuitionistic proof of Tarski's fixed point theorem making use of Separation, but no direct use of powerobjects is also presented in \cite{JM}. However, that proof appears to rely essentially on the assumption that $L$ is small, a requirement that by Theorem \ref{largestructures} is hardly met in systems without powersets.

To obtain a  version of the above corollary in a predicative system, we generalize the concept of bounded inductive definition from \cite{A86} to abstract inductive definitions. A \emph{bound} for an abstract inductive definition $\Phi$ is a set $\alpha$ such that, whenever $(b,a)\in \Phi$ there is $x\in \alpha$ such that the set $\down a$ is an image of $x$. An abstract inductive definition  $\Phi$ is \emph{bounded} if

\begin{enumerate}
  \item $\{b\in B\mid (b,a)\in \Phi\}$ is a set for every $a\in L$.
  \item $\Phi$ has a bound.
\end{enumerate}

Note that any abstract inductive definition that is a set is bounded. The following proposition generalizes the corresponding results for inductive definition on a set. Its proof is the obvious modification of the proof of \cite[Proposition 5.6]{AR}.  Exponentiation is used here for the first time in this paper.

\be{prop}[CZF$^-$+Exp]\label{boundedlocal}
Every bounded abstract inductive definition $\Phi$ is local.
\en{prop}

\be{rmk}\rm\label{}
If $\Phi$ is a set (rather than just bounded), a proof that it is local can in fact be given without any use of Exponentiation. Note that this in particular holds for standard inductive definitions.
\en{rmk}

\noindent Almost as obvious is the following generalization of  \cite[Proposition 5.3]{AR}, provable in CZF$^-$.

\be{prop}\label{JAsets}
If $\Phi$ is a local abstract inductive definition, $\texttt{$c_L$} J^{\alpha}$ and $\texttt{$c_L$} J^{\in \alpha}$ are sets for every $\alpha$.
\en{prop}

\be{proof}
We use induction on sets. Given $\alpha\in \mathrm{V}$, assume that for every $\beta\in \alpha$,  $\texttt{$c_L$} J^{\beta}$ is a set. Since, by Proposition \ref{c-operator},
$\texttt{$c_L$} J^{\in \alpha}\equiv \texttt{$c_L$} \bigcup_{\beta\in \alpha}J^{\beta} = \texttt{$c_L$} \bigcup_{\beta\in \alpha} \texttt{$c_L$} J^{\beta},$ and since $ \bigcup_{\beta\in \alpha} \texttt{$c_L$} J^{\beta}$ is a set, $\texttt{$c_L$} J^{\in \alpha}$ is a set, too.
Moreover, since for $Y\subseteq B$ a set,

\smallskip

\qquad \qquad \qquad \qquad
$\bar \Gamma_\Phi (Y)=\down \Gamma_\Phi (\bigvee Y),$
\smallskip

\noindent
 $\bar \Gamma_\Phi (Y)$ is a set for every set $Y$, as $\Phi$ is local and $L$ is set-generated. By Lemma \ref{classJ}, $\texttt{$c_L$} J^\alpha = \bar \Gamma_\Phi (\texttt{$c_L$} J^{\in \alpha})$, so that $\texttt{$c_L$} J^\alpha$ is  a set too. We can therefore conclude, by Set Induction, that for every $\alpha\in \mathrm{V}$, $\texttt{$c_L$} J^\alpha$, and so  $\texttt{$c_L$} J^{\in \alpha}$, are sets.
\en{proof}

In the case of a standard inductive definition $\Phi$, one shows that if $\Phi$ is bounded, then CZF + wREA proves that $I(\Phi)$ is a set \cite{AR}. We prove the corresponding result in our abstract context. Note that, in contrast with the  impredicative system CZF$^-$+Sep, CZF + wREA is a fully constructive theory \cite{A86}, the standard set-theoretic system for constructive predicative mathematics.

Recall that a  $\bigvee$-semilattice $L$ set-generated by a set $B$ is said to be \emph{set-presented} \cite{AR} if a mapping $D:B\to \Pow(\Pow(B))$ is given with the property that $b\leq \bigvee U\iff (\exists W\in D(b))W\subseteq U$, for every $b\in B, U\in \Pow(B)$.
The set-generated $\bigvee$-semilattice $\Pow(X)$, for $X$ a set, is for instance set-presented by the mapping $D(\{x\})=\{\{\{x\}\}\}$. The class of these $\bigvee$-semilattices is particularly  well-behaved (cf. \cite{A,AR} for more information).
In particular, provably in CZF+wREA, it includes a frame $L$ if and only if $L$ can be presented  by generators and relations
(in systems as those we are considering, even in the presence of full Sep, not every frame can be presented using generators and relations, cf. \cite{G1,CuriPA};  in particular, thus, Corollary \ref{TarskiTheorem} below will  not imply Corollary \ref{TarskiSeparation}).
In a system with powersets, every $\bigvee$-semilattice $L$ is trivially set-presented (by, for  $b\in B$, $D(b)=\{U\in\Pow(B)\mid b\leq \bigvee U\}$, with $B$ any base for $L$).

\be{thm}[CZF + wREA]\label{SetAbstractInductiveDef}
Let $\Phi$ be a bounded abstract inductive definition on a set-presented $\bigvee$-semilattice $L$.
Then, the smallest $\Phi-$closed class ${\cal I}(\Phi)$ is a set.
\en{thm}

\noindent To prove this result we shall need the following lemma.

\be{lem}\label{progressivity}
Given  an abstract inductive definition  $\Phi$ on a set-generated $\bigvee$-semilattice $L$, its associated iteration class $J$ satisfies, for every set $\alpha$,

\smallskip

\qquad \qquad $\texttt{$c_L$} J^{\in \alpha}\subseteq \texttt{$c_L$} J^{\alpha}.$
\en{lem}

\be{proof} It suffices to show that, for all $\alpha$,
 $ J^{\in \alpha}\subseteq c_L J^{\alpha}.$
We prove this by Set Induction. Let $x\in J^{\in \alpha}$, i.e., assume   $x\in J^\beta$ for $\beta\in \alpha$. Then, by Lemma \ref{classJ}, $x\in \bar \Gamma_\Phi(\texttt{$c_L$}J^{\in \beta})=\texttt{$c_L$} J^{\beta}$. Since, by inductive hypothesis,
$\texttt{$c_L$} J^{\in \beta}\subseteq \texttt{$c_L$} J^{ \beta}\subseteq \texttt{$c_L$}J^{\in \alpha},$ we conclude, by monotonicity of $\bar \Gamma_\Phi$ and Lemma \ref{classJ}, that $x\in \bar \Gamma_\Phi(\texttt{$c_L$}J^{\in \alpha})=\texttt{$c_L$} J^{\alpha}$.
\en{proof}

\be{proof}[Proof of Theorem \ref{SetAbstractInductiveDef}] Let $\alpha$ be a bound for $\Phi$, and let \smallskip

\qquad \qquad \qquad \qquad
$S=\alpha\cup\{V:(\exists b\in B) V\in D(b)\}= \alpha \cup \bigcup Range(D).$
\smallskip

\noindent
By Replacement (that is a consequence of Strong Collection) and Union, $S$ is a set. Then, by wREA, there is a weakly regular set $\alpha'$  such that $S\subseteq \alpha'$.

We claim that \smallskip

\qquad \qquad \qquad \qquad ${\cal I}(\Phi) \equiv \texttt{$c_L$} J^\infty = \texttt{$c_L$} J^{\in \alpha'},$\smallskip

\noindent the latter class being a set by Propositions \ref{boundedlocal}, \ref{JAsets}.

Since $J^{\in \alpha'}\subseteq J^\infty$, $\texttt{$c_L$} J^{\in \alpha'}\subseteq \texttt{$c_L$} J^\infty$.
So it remains to prove the converse, for which it suffices to show that $\texttt{$c_L$} J^{\in \alpha'}$ is $\Phi-$closed.
\smallskip

\noindent For $(b,a)\in \Phi$, assume $\down a\subseteq \texttt{$c_L$} J^{\in \alpha'}$. Since $\alpha$ is a bound for $\Phi$ there is a set $Z\in \alpha$ and an onto mapping $f:Z\twoheadrightarrow \down a$. So:
\smallskip

\qquad  $(\forall z\in Z)f(z)\in  \texttt{$c_L$} J^{\in \alpha'}$, i.e.,
\smallskip

\qquad $(\forall z\in Z)(\exists U\in \Pow(J^{\in \alpha'}))f(z)\leq \bigvee U$.
\smallskip

\noindent Since $L$ is set-presented, for every $z\in Z$, there is  $W\in D(f(z))$ such that $W\in \Pow(J^{\in \alpha'})$ and $f(z)\leq \bigvee W$. Therefore,
\smallskip

\qquad \qquad $(\forall d\in W)d\in  J^{\in \alpha'}$, i.e.,
\smallskip

\qquad \qquad $(\forall d\in W)(\exists \beta\in \alpha')d\in  J^{\beta}$,
\smallskip

\noindent which implies
\smallskip

\qquad \qquad $(\forall d\in W)(\exists \beta\in \alpha')d\in  \texttt{$c_L$} J^{\beta}$
\smallskip

\noindent (note that $\texttt{$c_L$} J^{\beta}$ is a set by Proposition \ref{JAsets}, while $J^{\beta}$ a priori need  not be). Now since, by construction, $W\in \alpha'$, and $\alpha'$ is weakly regular, there is a set $\gamma\in \alpha'$ such that
\smallskip

\qquad \qquad $(\forall d\in W)(\exists \beta\in \gamma)d\in  \texttt{$c_L$} J^{\beta}$.
\smallskip

\noindent Thus, $W\subseteq \bigcup_{\beta\in\gamma} \texttt{$c_L$} J^{\beta}$, so that $f(z)\in \texttt{$c_L$} \bigcup_{\beta\in\gamma} \texttt{$c_L$} J^{\beta}$. By Proposition \ref{c-operator}, $\texttt{$c_L$} \bigcup_{\beta\in\gamma} \texttt{$c_L$} J^{\beta}$ $ = \texttt{$c_L$} \bigcup_{\beta\in\gamma} J^{\beta} = \texttt{$c_L$}J^{\in \gamma}$. Further, by Lemma \ref{progressivity}, $\texttt{$c_L$}J^{\in \gamma}\subseteq \texttt{$c_L$}J^{ \gamma}$.
So we have shown
\smallskip

\qquad \qquad \qquad \qquad
$(\forall z\in Z)(\exists \gamma\in \alpha')f(z)\in  \texttt{$c_L$} J^{\gamma}.$

\smallskip

\noindent  Since $Z\in\alpha\subseteq \alpha'$, and  $\alpha'$ weakly regular, there is  $\delta\in\alpha'$ such that
\smallskip

\qquad \qquad \qquad \qquad
$(\forall z\in Z)(\exists \gamma\in \delta)f(z)\in  \texttt{$c_L$} J^{\gamma}.$
\smallskip

\noindent Therefore, $\down a\subseteq \bigcup_{\gamma\in\delta} \texttt{$c_L$} J^{\gamma}\subseteq  \texttt{$c_L$} J^{\in \delta}$.
By Lemma \ref{classJ}, then, $b\in \texttt{$c_L$}  J^{\delta}$. Finally, as $\delta\in \alpha'$, $\texttt{$c_L$}  J^{ \delta}\subseteq \texttt{$c_L$} J^{\in \alpha'}$, so that  $\texttt{$c_L$} J^{\in \alpha'}$ is $\Phi-$closed, as was to be proved.
\en{proof}

As a corollary we get the \emph{constructive Tarski's fixed point theorem}. Note that in classical set theory, or in a topos, every complete lattice is (set-generated and) set-presented, and  every monotone operator is obtained by a bounded abstract inductive definition (since in such a context any abstract inductive definition is a set).

\be{cor}[CZF + wREA]\label{TarskiTheorem} Let $\Gamma:L\to L$ be a monotone operator on
a set-presented $\bigvee$-semilattice $L$. If $\Gamma=\Gamma_\Phi$ for $\Phi$ be a bounded abstract inductive definition on $L$,
then $\Gamma$ has the least fixed point $p=\bigvee {\cal I}(\Phi)$.
\en{cor}

\section{Concluding remarks}
The  constructive derivation of Tarski's fixed point theorem we obtained may be perceived as strongly dependent on a particular formal system, the system CZF + wREA;  a more `system independent' derivation then might  seem desirable. In fact, one may regard  the presented proof in the following  perspective: the proofs of Theorem \ref{leastidclass} and of Theorem \ref{SetAbstractInductiveDef} show that the assertion that \emph{every bounded abstract inductive definition $\Phi$ on a set-presented $\bigvee$-semilattice $L$  inductively defines an element of $L$}  can be assumed as a constructively sound \emph{principle}, a generalization of
 the  result/principle that standard bounded inductive definitions on a set inductively define sets;
indeed, those proofs are derived in standard legitimate systems for constructive set theory (CZF$^-$ and its extension CZF + wREA). From this principle,  Tarski's theorem then immediately follows (Proposition \ref{idprefix}).

Abstract inductive definitions on a set-generated complete lattice generalize inductive definitions on a set (i.e., inductive definitions of the form $\Phi\subseteq X\times \Pow(X)$, for $X$ a set). Arbitrary inductive definitions can  be extended to an abstract setting in a similar way. The details of this extension will be presented in a forthcoming paper.

A version of Tarski's fixed point theorem for directed complete partial orders has been proved using either non-constructive arguments, or topos valid but impredicative reasoning (D. Pataraia). Constructive proofs of this and other fixed point theorems based on variants of the method presented in this paper will be the subject of future investigation.
Modifications of the method used in this paper could also allow for a solution of the problem, raised in \cite{A80} (see also \cite{CoqFixPoint}), of obtaining predicative constructions of Frege structures.

Finally: Tarski's theorem in its more general form states  that a monotone map on a complete lattice has a complete lattice of fixed points, so that in particular it has a greatest fixed point. A constructive version of this fact, based on the definition of an abstract (in the sense of this paper) counterpart of the concept of set/class  coinductively defined by an inductive definition, will also be investigated in future work.

\smallskip

\noindent \textbf{Acknowledgements}. I am grateful to the anonymous referees for their useful comments on a previous version of this paper.

\bibliographystyle{amsplain}

\end{document}